\documentclass[journal]{IEEEtran}
\ifCLASSINFOpdf
\else
\fi
\usepackage[linesnumbered,ruled]{algorithm2e}
\usepackage{amssymb}
\usepackage{nccmath}
\usepackage{cite}
\usepackage{graphicx}  
\usepackage[table]{xcolor}
\usepackage{verbatim}

\usepackage[utf8]{inputenc}
\usepackage{amsmath}
\usepackage{amsthm}     
\usepackage{accents}

\theoremstyle{remark}

\usepackage{tablefootnote}
\usepackage{url}

\DeclareMathOperator{\real}{Re}
\DeclareMathOperator{\imag}{Im}

\newcommand{\ubar}[1]{\underaccent{\bar}{#1}}
\begin{document}
%
\title{Network-Constrained Robust Unit Commitment for Hybrid AC/DC Transmission Grids}
%
%
%
\author{L.~P.~M.~I.~Sampath,
	M.~Hotz,~\IEEEmembership{Student~Member,~IEEE},
	H.~B.~Gooi,~\IEEEmembership{Senior~Member,~IEEE}
	and~W.~Utschick,~\IEEEmembership{Senior~Member,~IEEE}
	\thanks{L. P.~M.~I.~Sampath (Corresponding author) is with the Interdisciplinary Graduate School, Nanyang Technological University, Singapore 637371 (e-mail: mohashai001@e.ntu.edu.sg).} 
	\thanks{M.~Hotz and W.~Utschick are with the Department of Electrical and Computer Engineering, Technische Universit\"{a}t M\"{u}nchen, Munich D-80333, Germany (e-mail: matthias.hotz@tum.de; utschick@tum.de).}
	\thanks{H.~B.~Gooi is with the School of Electrical and Electronic Engineering, Nanyang Technological University, Singapore 639798 (e-mail: ehbgooi@ntu.edu.sg).}
	\thanks{This work is funded by the International Center of Energy Research	(ICER), established by Nanyang Technological University (NTU), Singapore	and Technische Universit\"{a}t M\"{u}nchen (TUM), Germany.}}

\maketitle

\begin{abstract}
The day-ahead energy and reserve management with transmission restrictions and voltage security limits is a challenging task for large-scale power systems in the presence of real-time variations caused by the uncertain demand and the fluctuating power output of renewable energy sources (RESs). The proposed formulation in this work supports joint scheduling of energy and reserve to promote an economic and reliable operation. 
To improve its scalability, a two-stage iterative optimization algorithm is proposed based on Bender's decomposition framework. Therewith, the optimal schedule is computed subject to the feasibility of AC network constraints (AC-NCs) at predetermined uncertain realizations. 
A convex relaxation is applied to AC-NCs to support the convergence of the algorithm. 
Moreover, the integration of RESs is often limited by transmission congestion issues in existing grids. For such grids, AC to DC conversion schemes can be viable and attractive options for capacity expansion. In this study, we adopt a hybrid AC/DC transmission grid (HTG) architecture with a specific topology to improve the effective utilization of the grid capacity and accommodate more demand and RESs, which is showcased via simulations on a large-scale network. 
Moreover, HTG induces exactness of the convex relaxation of AC-NCs, supporting the validity of the optimal schedule.
%
\end{abstract}

%
\IEEEpeerreviewmaketitle

\section{Introduction}\label{Intro}
\IEEEPARstart{U}{nit} commitment (UC) is a day-ahead optimal resource management problem to satisfy the forecasted demand subject to operational restrictions and the availability of generators~\cite{German2016,Yong2005,Nima2017,Amin2016}. Many countries have decided to integrate more renewable energy sources (RESs) into the generation mix due to environmental concerns~\cite{Garibaldi18}. The stochastic nature of the demand and RESs require power system operators to seek robust UC solutions against uncertainty sets which represent possible real-time variations in power injections.
Traditionally, reserves are allocated arbitrarily without analyzing the actual impact of uncertainties, which may lead to less security or less economic operations~\cite{German2016,Nima2017}. This motivates the joint planning of dispatch of energy and allocation of reserves within the UC formulation to schedule the most economic power supply mix that promotes an optimal operation of the system.

However, satisfying only the aggregated generation and demand balance may not be suffice in practice, as the line congestions may restrict the accessibility of geographically dispersed generation facilities to supply certain loads. Many notable studies utilize decoupled network constraints based on a simplified system model, which constitutes an easily tractable approximation of the network laws. 
However, the model mismatch requires more conservative system constraints in a practical implementation that leads to a suboptimal utilization~\cite{Castilo2016}. For instance, power losses can be significant and the active and reactive power may be strongly coupled, restricting the reactive power transmission which may lead to voltage security issues at some buses.
Consequently, the day-ahead schedule must be validated over the entire time horizon against the \emph{AC network constraints} (AC-NCs), which include nodal power balance, voltage bounds, power flow limits, and so forth~\cite{Hotz1,Yong2005,Nima2017,Amin2016}. 

The combined problem of UC with AC-NCs is extensive for practical transmission grids and typically cannot be solved directly using existing solvers~\cite{Amin2016,Yong2005}. For large-scale problems, Bender's decomposition framework (BDF) is a widely used iterative approach in the literature to improve scalability.
Many notable studies utilize BDF, with UC as a master problem (MP) and subproblems (SPs) to verify the feasibility of AC-NCs for possible real-time variations \cite{Yong2005,Nima2017,Amin2016,Azim2010}. 
In \cite{Yong2005}, UC is combined with AC-NCs for the forecasted scenario, ignoring the variability of RESs and demand. Recent studies propose multi-stage iterative approaches, using stochastic programs~\cite{Amin2016} and adaptive robust optimization~\cite{Nima2017} to address the network-constrained UC within the BDF under different uncertainty characterizations. However, the existing formulations either do not consider reserve allocations in inter-temporal constraints~\cite{Nima2017,Azim2010} or do not support the decomposition into SPs related to different time instants and uncertain realizations~\cite{Amin2016,Hongxing2016}. Therefore, the utilization of the former may not be viable in practice and the latter may not support the application to large-scale networks.%

In addition to the above, a proper formulation of AC-NCs invites nonconvexity into the optimization program due to the bilinear representation of voltage variables~\cite{Castilo2016}. 
The BDF requires convex SPs as the strong duality of SPs guarantees the validity of the Bender cuts and thereby the convergence of the algorithm~\cite{GBD}.
%
In that regard, several convex relaxations with improved mathematical and computational properties such as the semidefinite relaxation, second-order cone (SOC) relaxation and conic relaxation may be adopted in AC optimal power flow (OPF) problems, cf. 
\cite{LowConvex1, LowConvex2, Hotz1, Hotz2} and
the references therein. However, exactness of these relaxations is only guaranteed under certain conditions which are typically not fulfilled by traditional transmission grids~\cite{LavaeiConvex, LowConvex1, Hotz1}.

Transmission congestion can cause the dispatch of more expensive generators while cheaper generators are not fully utilized, due to the capacity limitation of some lines. According to the physical laws, the AC power flow (AC-PF) over a cycle of branches can stagnate due to congestion of a single branch in the cycle \cite{pablo}. The conventional approach to mitigate transmission congestion is the construction of new lines. However, this is often complicated and protracted by the construction of new corridors due to issues with right-of-way and public acceptance. 
As discussed in \cite{sanderACDC, AC_DC_ABB, ABB}, existing AC transmission lines can be converted to DC operation without corridor adjustments and the active power transmission capacity of the converted DC line can be two or more times the apparent power transmission capacity of the respective AC line, depending on the configuration. On this basis, a hybrid AC/DC transmission grid (HTG) architecture was recently proposed in \cite{Hotz1}, which serves as an upgrade strategy for traditional AC transmission grids to enhance its economic utilization and loadability \cite{Hotz1,Hotz2,Hotz3,MohashaPSCC}.
%
%
\subsection{Contributions and Outline}
This study aims to address the aforementioned issues by extending the authors' previous work in \cite{MohashaPSCC} to day-ahead energy and reserve scheduling on large-scale power systems under uncertainty using the BDF. 
In contrast to the state-of-art techniques~\cite{Yong2005,Nima2017,Amin2016,Azim2010} that combine UC with AC-NCs, the main contributions of this study are as follows.


\begin{enumerate}
\item Contributions to the problem formulation:
The MP formulation of commitment and generation scheduling of units is extended with the scheduling of reserves, which comply with inter-temporal ramping and transmission capacity restrictions while offering robust operation within the confidence bounds of uncertain generation and demand. 
\item Contributions to the solution method:
%
%
Decomposition of SPs with respect to the sets of scenarios and time instants using the BDF, enabling a parallel computation of SPs. 
In addition, SPs are formulated in an effective and unified manner. 
Therein, a low number of feasibility slack variables is used to construct tight Bender's feasibility cuts (BFCs) and compact SPs that reduce the memory requirements. Further, SPs are cast as standard OPF problems promoting the utilization of existing OPF packages.
Moreover, an SOC relaxation is applied to AC-NCs to form \emph{convex} SPs which validate the BFCs and ensure convergence of the iterative solution method. 
\item Contributions in the results:
In the simulations, a large-scale network with RES integration is considered. To the best of the authors' knowledge, this is lacking among the studies that combine UC with AC-NCs in the literature. Case studies showcase the effectiveness of a topology preserving capacity expansion procedure \cite{Hotz2} and measure its merits in accommodating the growing demand and improving the hosting capacity for RESs. 
\end{enumerate}
%
%
The paper is organized as follows. The day-ahead optimal energy and reserve management problem is formulated in Section~\ref{RO}. The decomposition of the combined UC with AC-NCs is explained in Section~\ref{GBD}. Section~\ref{TwoStgAlgo} discusses the convex relaxation of SPs and introduces a two-stage iterative algorithm to compute the optimal day-ahead schedule. 
Simulations results are presented in Section~\ref{Results}, followed by the conclusion in Section~\ref{Conclu}.
\section{Robust Optimization Model for Day-Ahead Market Planning}\label{RO}
The day-ahead market planning involves the preparation of the system to operate with a minimum generation cost schedule while assuring the operational feasibility for worst-case real-time operating conditions. 
In this section, the formulation of a scenario-based robust optimization model is presented below.
In general, the day-ahead optimal planning involves decisions over a predefined time horizon $ \mathcal{T}=\{1,\ldots,T\} $.
\subsection{Uncertainty Characterization} \label{TOAT}
The RES power outputs and the load demand are uncertain in advance. 
There are several forecasting techniques which employ historical data to model and predict these stochastic and time-varying outcomes. In this paper, the uncertain variables are assumed to be varying within a polyhedron. This polyhedron is defined by the deterministic confidence bounds for each uncertain variable over the time horizon which can be computed based on historical data using statistical inference techniques~\cite{Pinson2010}. 
The feasibility of the solution against the worst-case scenario can be found by maximizing the constraint violations over the polyhedral uncertainty set~\cite{Hongxing2016}. Considering uncertainties pertained to bus injections, the worst-case scenario should be among the extreme points of the uncertainty polytope \cite[Sec. III-C]{YZhang2013}. However, exploring all extreme points is computationally expensive as their number rises exponentially with the dimension of the uncertainty set~\cite{Xiang16}. In this context, \emph{Taguchi's orthogonal array testing} (TOAT) is a widely used special set of \emph{orthogonal arrays} (OAs) which provides a fractional amount of the full set of extreme points of the polytope~\cite{Wang2017,HanYu2012,TOAT,Xiang16}. The reader may refer to~\cite{TOAT} for an extensive review on the application of the TOAT method for power system uncertainty characterization.%

Let $ \mathcal{N} $ be the set of buses in the system, where $ |\mathcal{N}|=N $, and let $ \Omega $ be the set of uncertain scenarios provided by the TOAT method. Furthermore, let $ P^{t,\omega}_{{\rm D},n} $ be a realization of the demand for scenario~$ \omega \in \Omega$ at bus~$ n \in \mathcal{N} $ at time~$ t \in \mathcal{T} $ that corresponds to an extreme point of the polyhedral uncertainty set, \emph{i.e.}, a point on the boundary of the confidence interval.%
\begin{equation}\label{Pd}%
P^{t,\omega}_{{\rm D},n} \in \Big\{{P}^{t}_{{\rm D},n} - P^{t,\zeta-}_{{\rm D},n},\ {P}^{t}_{{\rm D},n} +P^{t,\zeta+}_{{\rm D},n} \Big\};\, \forall \omega \in \Omega \backslash \{0\}%
\end{equation}%

Therein, the nonnegative parameters $ {P}^{t}_{{\rm D},n},\ P^{t,\zeta-}_{{\rm D},n}$ and $ P^{t,\zeta+}_{{\rm D},n} $ are the mean and the expected deviations for the upper and lower confidence bounds of the active power demand at bus~$ n $ at time $ t $, respectively. Similarly, the extreme points for the RES outputs $ P^{t,\omega}_{{\rm W},n} $ can be defined using the mean value $ {P}^{t}_{{\rm W},n} $ and the deviations $ P^{t,\zeta+}_{{\rm W},n}$ and $P^{t,\zeta-}_{{\rm W},n} $. Here, $ \omega=0 $ represents the base-case scenario at which demand and RES outputs assume their forecast values, \emph{i.e.}, $ {P}^{t,0}_{{\rm D},n}={P}^{t}_{{\rm D},n} $ and $ {P}^{t,0}_{{\rm W},n}={P}^{t}_{{\rm W},n} $.%
\subsection{Unit Commitment and Power Dispatch Formulation}
The proposed day-ahead energy and reserve management problem in~\eqref{UC_Prob} derives the optimal schedule for the \emph{base-case} operation. 
The optimal solution include the generator commitment plan, power dispatch values, and the reserve margins for all the generators in the system. For simplicity of exposition, we assume one generator per bus in the system. The decision variables of UC pertaining to generator~$ n \in \mathcal{N} $ at time~$ t \in \mathcal{T} $ are defined as follows.
\begin{equation}\label{VAR}
\Phi^t_n=\Big\{ \hat{u}_n^t, \check{u}_n^t, u_n^t, P_{\mathrm{G},n}^{t}, r^{t}_{{\rm u},n}, r^{t}_{{\rm d},n} \Big\};\,    \forall n\in{\mathcal{N}},\, \forall t\in{\mathcal{T}}
\end{equation}
Therein, $ \hat{u}_n^t$ and $ \check{u}_n^t $ are the startup and shutdown instances, respectively, $ u_n^t$ is the commitment state, $ P_{\mathrm{G},n}^{t} $ is the power dispatch, and $ r^{t}_{{\rm u},n}$ and $ r^{t}_{{\rm d},n} $ denote the up and down reserve allocations, respectively.
\begin{subequations}\label{UC_Prob}
\begin{align}
\min_{\substack{\Phi^\mathcal{T}_\mathcal{N}}}\ & 
\sum_{t\in{\mathcal{T}}} \sum_{n\in{\mathcal{N}} } \left[ c^{\rm e}_{1,n} P_{\mathrm{G},n}^{t}+c^{\rm e}_{0,n}{u}_n^t + c^{\rm s}_n  \hat{u}_n^t + c^{\rm d}_n  \check{u}_n^t \right] \label{UC_obj} \\
{\rm s.t.}\ & u_n^t,\ \hat{u}_n^t,\ \check{u}_n^t \in \{0,1\};\,    \forall n\in{\mathcal{N}},\,  \forall t\in{\mathcal{T}}\\
& \hat{u}_n^t - \check{u}_n^t = u_n^t - u_n^{t-1} ;\, \forall n\in{\mathcal{N}},\, \forall t\in{\mathcal{T}}\hspace{0.2mm} \backslash \{1\} \label{susd}\\
&\sum_{\tau =t}^{t+\hat{T}_n-1} u_{n}^{\tau} \geq \hat{T}_n \hat{u}_n^t;\,  \forall n \in \mathcal{N}, \, \forall t \in \mathcal{T}  \label{UT} \\
&\sum_{\tau=t}^{t+\check{T}_n-1} \left[ 1 - u_{n}^{{\tau}}\right]  \geq \check{T}_n \check{u}_n^t;\,  \forall n \in \mathcal{N}, \, \forall t \in \mathcal{T} \label{DT} \\
& P_{{\rm G},n}^{t}+ r^{t}_{{\rm u},n} \leq u_n^t P_{{\rm G},n}^{\max}  ;\, \forall n\in{\mathcal{N}},\, \forall t\in{\mathcal{T}} \label{Gub}\\
& P_{{\rm G},n}^{t} - r^{t}_{{\rm d},n} \geq u_n^t P_{{\rm G},n}^{\min} ;\,    \forall n\in{\mathcal{N}},\,  \forall t\in{\mathcal{T}} \label{Glb}\\
& 0 \leq r^{t}_{{\rm u},n} \leq u_n^{t} R^\delta_{{\rm u},n}  ;\,    \forall n\in{\mathcal{N}},\, \forall t\in{\mathcal{T}} \label{resub}\\
& 0 \leq r^{t}_{{\rm d},n} \leq u_n^{t} R^\delta_{{\rm d},n}  ;\,    \forall n\in{\mathcal{N}},\,  \forall t\in{\mathcal{T}} \label{reslb}\\
& P_{{\rm G},n}^{t} + r^{t}_{{\rm u},n} - P_{{\rm G},n}^{t-1} + r^{t-1}_{{\rm d},n} \leq u_n^{t-1}R^\Delta_{{\rm u},n}  \nonumber \\
& \hspace{0.75cm} + \left( 1-u_n^{t-1}\right) P_{{\rm G},n}^{\min};\, \forall n \in \mathcal{N}, \, \forall t \in \mathcal{T}\hspace{0.2mm} \backslash \{1\} \label{Rup}\\ 
& P_{{\rm G},n}^{t-1}+ r^{t-1}_{{\rm u},n} - P_{{\rm G},n}^{t} + r^{t}_{{\rm d},n} \leq u_n^{t}R^\Delta_{{\rm d},n}  \nonumber\\
& \hspace{0.75cm} + \left( 1-u_n^{t}\right) P_{{\rm G},n}^{\min} ;\, \forall n \in \mathcal{N}, \, \forall t \in \mathcal{T}\hspace{0.2mm} \backslash \{1\} \label{Rdn}\\
%
%
& \sum_{n \in \mathcal{N}}P_{{\rm G},n}^{t} \geq \sum_{n \in \mathcal{N}}P_{{\rm D},n}^{t} - \sum_{n \in \mathcal{N}}P_{{\rm W},n}^{t}+P^{t}_{\rm loss} ;\, \forall t \in \mathcal{T} \label{Ebal}\\ 
& \sum_{n \in \mathcal{N}}r^{t}_{{\rm u},n} \geq \sum_{n \in \mathcal{N}} \Big[  \max_{\omega \in {\Omega}} P_{{\rm D},n}^{t,\omega}-P_{{\rm D},n}^{t}\Big] +\max_{\omega \in {\Omega}}\Delta P^{t,\omega}_{\rm loss} \nonumber \\
&\hspace{1.5cm} +\sum_{n \in \mathcal{N}}\Big[ P_{{\rm W},n}^{t} - \min_{\omega \in {\Omega}} P_{{\rm W},n}^{t,\omega} \Big] ;\, \forall t \in \mathcal{T} \label{upRes}\\
& \frac{1}{\alpha^t} \sum_{n \in \mathcal{N}}r^{t}_{{\rm d},n} \geq \sum_{n \in \mathcal{N}} \Big[ P_{{\rm D},n}^{t}- \min_{\omega \in {\Omega}} P_{{\rm D},n}^{t,\omega}\Big] -\min_{\omega \in {\Omega}}\Delta {P}^{t,\omega}_{\rm loss} \nonumber \\
&\hspace{1.5cm} +\sum_{n \in \mathcal{N}}\Big[ \max_{\omega \in {\Omega}} P_{{\rm W},n}^{t,\omega} - P_{{\rm W},n}^{t} \Big] ;\, \forall t \in \mathcal{T} \label{dnRes}
\end{align}
\end{subequations}
In \eqref{UC_Prob}, $ \Phi^\mathcal{T}_\mathcal{N} $ denotes $ N|\mathcal{T}| $-tuple $ \Phi^\mathcal{T}_\mathcal{N}=(\Phi^t_n)^{t \in \mathcal{T}}_{n \in \mathcal{N}} $. Further, $ \Delta P^{t,\omega}_{\rm loss}=P^{t,\omega}_{\rm loss}-P^{t}_{\rm loss};\, \forall \omega \in {\Omega},\, \forall t \in \mathcal{T} $; $ P^{t}_{\rm loss} $ and $ P^{t,\omega}_{\rm loss} $ are the total power losses at the base-case and at scenario~$ \omega $, which can be calculated using~\eqref{P_LOSS}.
The objective function of the UC in~\eqref{UC_obj} is the minimization of the generation cost over the specified time horizon $ \mathcal{T} $. The generation cost function comprises the marginal and fixed energy cost terms $ c^{\rm e}_{1,n} $ and $ c^{\rm e}_{0,n} $ as well as the startup cost term $ c^{\rm s}_{n} $ and shutdown cost term $ c^{\rm d}_{n} $. Equation~\eqref{susd} explains the relation of the binary variables $ \hat{u}_n^t$, $ \check{u}_n^t$ and $ u_n^t $. Minimum uptime $ \hat{T}_n $ and downtime $ \check{T}_n $ requirements of all the generation units are enforced with \eqref{UT} and \eqref{DT} respectively. 
The energy and reserve allocations are constrained by \eqref{Gub} and \eqref{Glb}, respecting the generation capability bounds.
Similarly, \eqref{resub} and \eqref{reslb} limit the up and down reserve capability within the short-term up and down ramp rates $ R^\delta_{{\rm u},n} $ and $ R^\delta_{{\rm d},n} $ respectively. Moreover, reserve allocations in between time instances are consistent with the inter-temporal startup and shutdown as well as the operating ramp rates $ R^\Delta_{{\rm u},n} $ and $ R^\Delta_{{\rm d},n} $ as per \eqref{Rup} and \eqref{Rdn} respectively. This is a salient feature in the proposed formulation compared to the existing studies~\cite{Nima2017,Amin2016,Azim2010}. 
The minimum generation requirement for the base-case operation is satisfied by~\eqref{Ebal} over the time frame. 
Similarly, the minimum required reserve allocation can be determined by \eqref{upRes} and \eqref{dnRes}, considering the possible fluctuations in demand and RES outputs. 
The up-reserve allocations typically increase the base-case operating cost as it usually causes an inferior utilization of economical units for the base-case energy dispatch.  
Moreover, $ \alpha^t $ is a constant parameter with the default value of $1$ for all $ t \in \mathcal{T} $. Its relevance to the solution process is explained in Section~\ref{minPL}.%
\subsection{AC Network Constraint Formulation}
The day-ahead generator commitment and power dispatch schedule in~\eqref{UC_Prob} is subjected to the real-time variations in the demand and RES outputs. 
%
Let $ \mathcal{K}$ and $ \mathcal{L}$ be the set of AC lines and the set of DC lines of the power system respectively, where $ |\mathcal{K}| = K$ and $ \  |\mathcal{L}| = L  $. 
Therewith, the actual active and reactive power dispatch values $ P^{t,\omega}_{\mathrm{G},n} $ and $ Q^{t,\omega}_{\mathrm{G},n} $, the complex bus voltage vector $ v^{t,\omega} \in \mathbb{C}^{N} $  and DC branch flow vector $ p^{t,\omega} \in \mathbb{R}^{L} $ can be used to formulate the AC-NCs for scenario~$ \omega \in \Omega$  at time~$ t \in \mathcal{T}$ as in~\eqref{ACPF1}. The reader may refer to \cite[Sec.~II and Sec.~III]{Hotz1} for the derivation of the constraints from first principles.
\begin{subequations}\label{ACPF1}
	\begin{align}
	& P^{t,\omega}_{\mathrm{G},n} - P^{t,\omega}_{\mathrm{D},n} +P^{t,\omega}_{\mathrm{W},n} \nonumber \\
	&\hspace{1.5cm} = ({v}^{t,\omega})^{\mathrm{H}} {P}_n {v}^{t,\omega} + {h}_n^{\mathrm{T}} {p}^{t,\omega} ;\   \forall n\in{\mathcal{N}} \label{P_bal} \\ 
	& Q^{t,\omega}_{\mathrm{G},n}- Q^{t,\omega}_{\mathrm{D},n} +Q^{t,\omega}_{\mathrm{W},n} = (v^{t,\omega})^{\mathrm{H}} {Q}_n {v}^{t,\omega};\   \forall n\in{\mathcal{N}} \label{Q_bal} \\
	& (v_n^{\min})^2 \leq (v^{t,\omega})^{\mathrm{H}} {M}_n v^{t,\omega} \leq (v_n^{\max})^2;\  \forall n\in{\mathcal{N}} \label{V_lim}\\
	& (v^{t,\omega})^{\mathrm{H}} {\hat{I}}_{k} v^{t,\omega} \leq (I_{k}^{\max})^2 ;\   \forall k \in{\mathcal{K}} \label{I_fl1} \\
	& (v^{t,\omega})^{\mathrm{H}} {\check{I}}_{k} v^{t,\omega} \leq (I_{k}^{\max})^2 ; \  \forall k \in{\mathcal{K}} \label{I_fl2} \\
	& p_l^{\min} \leq p^{t,\omega}_{l} \leq p_l^{\max} ;\   \forall l \in\mathcal{L} \label{DClim1} \\
	& -{r}^{t,\omega}_{{\rm d},n} \leq P^{t,\omega}_{\mathrm{G},n}-{P}^{t}_{\mathrm{G},n} \leq {r}^{t,\omega}_{{\rm u},n} ;\  \forall n\in{\mathcal{N}} \label{RES_lim1}\\
	& {u}^t_n Q^{\min}_{\mathrm{G},n} \leq Q^{t,\omega}_{\mathrm{G},n} \leq {u}^t_n Q^{\max}_{\mathrm{G},n} ;\  \forall n\in{\mathcal{N}} \label{Qlim1} 
	\end{align}
\end{subequations}

In \eqref{ACPF1}, $(\cdot)^{\mathrm{T}}$ and $(\cdot)^{\mathrm{H}}$ denote the transpose and Hermitian transpose, respectively. 
Constraints \eqref{P_bal} and \eqref{Q_bal} constitute the active and reactive power balance equation at bus $ n $ respectively. The matrices $ {P}_n$ and $ {Q}_n \in \mathbb{S}^{N} $ are functions of the bus admittance matrix for bus~$ n $, while $ {h}_n \in \mathbb{R}^{L} $ describes the power flow, losses (assumed to be proportional to the power flow through DC lines), and connectivity of DC lines to bus~$ n $. Hence, the first and second term of \eqref{P_bal} define AC and DC extractions of bus~$ n $, respectively. Constraint~\eqref{V_lim} explains the bus voltage magnitude limits at bus $ n $ characterized by the matrix $ {M}_n $. Constraints \eqref{I_fl1} and \eqref{I_fl2} limit the bidirectional current flow through AC line $ k $, where $ \hat{I}_{k}, \ \check{I}_{k} \in \mathbb{S}^N $ characterize the respective functions of the branch admittance matrix. Power flow limits of DC line~$ l $ are incorporated in~\eqref{DClim1}. The active and reactive power generation capability at bus~$ n $ is represented by \eqref{RES_lim1} and \eqref{Qlim1}, respectively. Up and down reserve allocations $ {r}^{t,\omega}_{{\rm u},n},\ {r}^{t,\omega}_{{\rm d},n} $ for scenario $ \omega $ at time $ t $ are assigned as below. 
\begin{equation} \label{r_u}
\big\{{r}^{t,\omega}_{{\rm u},n},\,{r}^{t,\omega}_{{\rm d},n}\big\} = 
\left\{
\begin{aligned}%
& \big\{{r}^{t}_{{\rm u},n},\,{r}^{t}_{{\rm d},n}\big\}           ; \quad \omega \in \Omega \backslash \{0\}  \\
&\hspace{0.45cm} \{0,\,0\} \hspace{0.45cm} ; \quad \omega=0 \hspace{7.5mm} 
\end{aligned}%
\right. 
\end{equation}%
The total power loss can be computed as in \eqref{P_LOSS}.%
\begin{equation}\label{P_LOSS}%
P^{t,\omega}_{\rm loss}=\sum_{n \in \mathcal{N}}\Big[ P_{{\rm G},n}^{t,\omega} - P_{{\rm D},n}^{t,\omega} + P_{{\rm W},n}^{t,\omega} \Big]%
\end{equation}%
%
\vspace{-0.5em}%
\section{Problem Decomposition}\label{GBD}
It is evident from Section~\ref{RO} that the network-constrained robust day-ahead energy and reserve management problem is a combination of \eqref{UC_Prob} and the AC-NCs \eqref{ACPF1}, $ \forall \omega \in \Omega $, $ \forall t \in \mathcal{T} $. The dimensionality of the combined problem is proportional to $ N|\mathcal{T}||\Omega| $ which is extensive for a direct solution approach. To this end, a two-stage iterative approach is proposed in this section based on \emph{Bender's decomposition}. The latter is a widely used and effective framework for decomposing large-scale mixed-integer optimization problems with guaranteed $ \epsilon $-convergence in finite number of iterations, provided that given conditions, e.g.~\cite[Theorem~2.5]{GBD}, hold. %
In general, the BDF comprises a master problem (MP) and several convex subproblems (SPs) which sequentially iterate in the solution process. The formulation in Section~\ref{RO} shows that the structure of AC-NCs is consistent for all the scenarios~$ \omega \in \Omega $ and time instances~$ t \in \mathcal{T} $. 
Therewith, the MP and SPs are formulated as follows.
\vspace{-0.5em}%
\subsection{Master Problem Formulation}
In the first stage, the day-ahead energy and reserve scheduling problem~\eqref{UC_Prob} at fixed $ P^{t,\omega}_{\rm loss}, \, {\forall \omega \in \Omega}, \, \forall t \in \mathcal{T} $, is considered as the MP, which computes the day-ahead decisions $ \Phi^\mathcal{T}_\mathcal{N} $.
Here, the MP formulation accounts for all the inter-temporal constraints and the corresponding variables which include the generator commitment decisions, the base-case power dispatch and the reserve allocations, enabling a full decomposition of AC-NCs~\eqref{ACPF1} per time $ t $ and scenario $ \omega $.%
\vspace{-0.5em}%
\subsection{Subproblem Formulation}
The second stage involves SPs based on AC-NCs \eqref{ACPF1}, $ \forall \omega \in \Omega $, $ \forall t \in \mathcal{T} $, in order to verify the worst-case operational feasibility of the day-ahead decisions. Hereafter, the day-ahead decisions $ \big\{ u_n^t,\, P^{t}_{\mathrm{G},n},\, r^{t}_{{\rm u},n},\, r^{t}_{{\rm d},n} \big\},\,    \forall n\in{\mathcal{N}},\, \forall t\in{\mathcal{T}} $ of \eqref{UC_Prob} which interact with the AC-NCs are referred as \emph{complicating variables}. These are fixed in the coupling constraints \eqref{RES_lim1} and \eqref{Qlim1} at the given values from the MP. 
%
Therefore, a particular day-ahead schedule $ {\Phi}^\mathcal{T}_\mathcal{N} $ may induce violations in voltage bounds and/or power flow limits, and/or the active power dispatch and/or reserves may be insufficient. 
Hence, some constraints must be relaxed using slack variables and the use of slack needs to be penalized in the objective to find the solution with minimal constraint violations \cite{Amin2016,Nima2017}. In this study, the real-time feasibility verification SP in~\eqref{SubProb} for scenario~$ \omega \in \Omega $ at time~$ t \in \mathcal{T} $ is formulated using AC-NCs with flexible injection bounds, \emph{i.e.}, the coupling constraints \eqref{RES_lim1} and \eqref{Qlim1} are relaxed. 

To this end, let a convex piece-wise linear function ${ \sigma^t_n : \mathbb{R}^3 \rightarrow \mathbb{R}} $ comprising three linear segments with positive, zero and negative gradients be defined as
\begin{equation} \label{SIGMA_x}%
\sigma^{t}_n\hspace{-1mm}\left( x^{t}_n,\bar{x}^{t}_n,\ubar{x}^{t}_n\right)\hspace{-1mm} = \hspace{-1mm}
\left\{
\begin{aligned}\hspace{0mm}%
& \gamma(x^{t}_n -  \bar{x}^{t}_n)         \hspace{4.3mm}   \,  ; \,   \bar{x}^{t}_n < & \hspace{-2.5mm} {x}^{t}_n &  \\
&\hspace{0.9cm} 0       \hspace{1.05cm}   \,    ; \,  \ubar{x}^{t}_n    \leq  & \hspace{-2.5mm}  {x}^{t}_n &  \leq \bar{x}^{t}_n \\
& -\gamma({x}^{t}_n - \ubar{x}^{t}_n) \,   ; \, & \hspace{-2.5mm} {x}^{t}_n &   < \ubar{x}^{t}_n    
\end{aligned}%
\right. 
\end{equation}%
where $ \gamma > \underset{n \in \mathcal{N}}{\max} \, c^e_{1,n} $. Therewith, the proposed SP reads
%
%
\begin{subequations}\label{SubProb}
	\begin{align}
	z^{t,\omega}=\min_{\substack{P^{t,\omega}_{\mathrm{G},n},\, Q^{t,\omega}_{\mathrm{G},n}\\ {v}^{t,\omega},\,{p}^{t,\omega}}}\, & 
	\sum_{n\in{\mathcal{N}}} \sigma^{t}_n\hspace{-1mm}\left( P^{t,\omega}_{\mathrm{G},n},\, {P}^{t}_{\mathrm{G},n} +{r}^{t,\omega}_{{\rm u},n}, \, {P}^{t}_{\mathrm{G},n} -{r}^{t,\omega}_{{\rm d},n}\right)  \nonumber \\
	& \hspace{-3.5mm}+ \sum_{n\in{\mathcal{N}}} \sigma^{t}_n\hspace{-1mm}\left( Q^{t,\omega}_{\mathrm{G},n},\,{u}^t_n Q^{\rm max}_{\mathrm{G},n},\, {u}^t_nQ^{\rm min}_{\mathrm{G},n}\right) \label{SP_obj} \\
	{\rm s.t.}\ & \eqref{P_bal}{\text -}\eqref{DClim1}
	\end{align}
\end{subequations}%

The bounds on $ P^{t,\omega}_{\mathrm{G},n} $ and $ Q^{t,\omega}_{\mathrm{G},n} $ in \eqref{RES_lim1} and \eqref{Qlim1} are relaxed using soft limits in the objective function~\eqref{SP_obj}. In contrast to the existing literature \cite{Amin2016,Nima2017}, this formulation of SPs requires neither constraint manipulation nor additional variables, \emph{i.e.}, SPs are cast as standard OPFs which enables the utilization of existing OPF packages. 
The SP in~\eqref{SubProb} computes the minimum amount of constraint violations $ z^{t,\omega} $ necessitated by the day-ahead schedule for scenario~$ \omega $ at time~$ t $. The corresponding $ P^{t,\omega}_{\rm loss} $ can be calculated for each SP using~\eqref{P_LOSS}.
\subsection{Bender Cut Formulation}\label{Sec_BenCut}
In BDF, the dual domain of SPs is utilized to deduce the sensitivities of complicating variables on constraint violations. Therewith, Bender's feasibility cuts (BFCs) are constructed in every iteration and integrated into the MP for subsequent iterations. Here, the role of BFCs is to rectify $ \Phi^{\mathcal{T}}_{\mathcal{N}} $ in order to avoid constraint violations in AC-NCs.
%

Since \eqref{RES_lim1} and \eqref{Qlim1} are omitted in the SP~\eqref{SubProb}, the corresponding dual variables of reserve allocations and generator commitment are extracted as follows. 
Let $ \pi^{t,\omega}_{{\rm P},n} $ be the dual variable of active power balance constraint \eqref{P_bal} and $ \hat{\pi}^{t,\omega}_{{\rm R},n} $ and $ \check{\pi}^{t,\omega}_{{\rm R},n} $ be the dual variables of the up- and down-reserve capability bounds \eqref{RES_lim1} of bus~$ n \in \mathcal{N}$ for scenario~$\omega \in \Omega \backslash \{0\}$ at time~$ t \in \mathcal{T} $ respectively. Then, using the tolerance parameter $ \epsilon > 0 $ for bound satisfaction, we identify that
\begin{itemize}
	\item If $ P^{t,\omega}_{\mathrm{G},n} > {P}^{t}_{\mathrm{G},n} + {r}^{t}_{{\rm u},n}+\epsilon $, then
	$ \hat{\pi}^{t,\omega}_{{\rm R},n}={\pi}^{t,\omega}_{{\rm P},n} $,
	$ \check{\pi}^{t,\omega}_{{\rm R},n}=0 $,
	\item else if $ P^{t,\omega}_{\mathrm{G},n} < {P}^{t}_{\mathrm{G},n} -{r}^{t}_{{\rm d},n} -\epsilon $, then
	$ \hat{\pi}^{t,\omega}_{{\rm R},n}=0 $, 
	$ \check{\pi}^{t,\omega}_{{\rm R},n} =-{\pi}^{t,\omega}_{{\rm P},n} $,
	\item otherwise
	$ \hat{\pi}^{t,\omega}_{{\rm R},n} =\check{\pi}^{t,\omega}_{{\rm R},n} = 0 $.
\end{itemize}
%
%
%
%
Similarly, let $ \pi^{t,\omega}_{{\rm Q},n} $ be the dual variable of reactive power balance constraint \eqref{Q_bal} and $ \hat{\pi}^{t,\omega}_{{\rm Q},n} $ and $ \check{\pi}^{t,\omega}_{{\rm Q},n} $ be the dual variables of reactive power capability upper and lower bounds \eqref{Qlim1} of bus~$ n \in \mathcal{N}$ for scenario~$\omega \in \Omega$ at time~$ t \in \mathcal{T} $ respectively.
\begin{itemize}
	\item If $ Q^{t,\omega}_{\mathrm{G},n} > {u}^t_n Q^{\rm max}_{\mathrm{G},n} +\epsilon $, then
	$ \hat{\pi}^{t,\omega}_{{\rm Q},n}={\pi}^{t,\omega}_{{\rm Q},n} $,
	$ \check{\pi}^{t,\omega}_{{\rm Q},n}=0 $
	\item else if $ Q^{t,\omega}_{\mathrm{G},n} < {u}^t_nQ^{\rm min}_{\mathrm{G},n} -\epsilon $, then
	$ \hat{\pi}^{t,\omega}_{{\rm Q},n}=0 $, 
	$ \check{\pi}^{t,\omega}_{{\rm Q},n} =-{\pi}^{t,\omega}_{{\rm Q},n} $
	\item otherwise
	$ \hat{\pi}^{t,\omega}_{{\rm Q},n} =\check{\pi}^{t,\omega}_{{\rm Q},n} =0 $.
\end{itemize}
%

Then, the dual variables and $ z^{t,\omega} $ for all scenarios~$ \omega \in \Omega $ are averaged at each time~$ t $ as in~\eqref{AvgFun}, in order to reduce the number of constraints imposed on MP at every iteration~\cite{Amin2016}.
\begin{equation}\label{AvgFun}
{z}^{t}= z^{t,0}+\frac{1}{|\Omega|-1}\sum_{\omega \in \Omega \setminus \{0\}}z^{t,\omega}\, ; \ \forall t \in \mathcal{T}
\end{equation}
Similarly, $ {\pi}^{t}_{{\rm P},n}, \, \hat{\pi}^{t}_{{\rm R},n}, \, \check{\pi}^{t}_{{\rm R},n}, \, \hat{\pi}^{t}_{{\rm Q},n} $ and $ \check{\pi}^{t}_{{\rm Q},n} $ can be computed. %
To this end, BFCs are formulated as in~\eqref{BenCut}.
\begin{align}\label{BenCut}
{z}^t & +\hspace{-1mm}\sum_{n\in{\mathcal{N}}}\hspace{-0.5mm}\Big[\hat{\pi}^{t}_{{\rm Q},n} Q^{\max}_{\mathrm{G},n}-\check{\pi}^{t}_{{\rm Q},n} Q^{\min}_{\mathrm{G},n}\Big]\Big({u}^{t,(\eta)}_{n}- {u}^{t,(\eta-1)}_{n}\Big) \nonumber\\ 
& +\hspace{-1mm}\sum_{n\in{\mathcal{N}}}\hspace{-0.5mm}\hat{\pi}^{t}_{{\rm R},n} \Big( {r}^{t,(\eta)}_{{\rm u},n}-{r}^{t,(\eta-1)}_{{\rm u},n} \Big) 
+\hspace{-1mm}\sum_{n\in{\mathcal{N}}}\hspace{-0.5mm}\check{\pi}^{t}_{{\rm R},n} \Big( {r}^{t,(\eta)}_{{\rm d},n} -{r}^{t,(\eta-1)}_{{\rm d},n} \Big) \nonumber\\
& +\hspace{-1mm}\sum_{n\in{\mathcal{N}}}\hspace{-0.5mm}{\pi}^{t}_{{\rm P},n} \Big({P}^{t,(\eta)}_{\mathrm{G},n} -{P}^{t,(\eta-1)}_{\mathrm{G},n}\Big) \leq 0 \, ; \ \forall t \in \mathcal{T}
\end{align}
Therein, 
$ \eta $ is the current iteration of the algorithm. 
The current variables pertaining to the MP in~\eqref{UC_Prob} are termed with the superscript~$ \eta $. The BFCs are included in the MP to rectify the commitment, generation and reserve scheduling to eliminate violations in AC-NCs.
\section{Two-stage Iterative Algorithm based on Bender's Decomposition}\label{TwoStgAlgo}
\subsection{Convex Relaxation of the Subproblems} \label{ConvRel}
The SP in~\eqref{SubProb} are nonconvex. As explained in Section~\ref{Sec_BenCut}, the sensitivities of the complicating variables are used to feedback the information on constraint violations to the MP in the form of BFCs. Here, the sensitivities are computed based on the dual variables of the SPs. In that respect, if the SPs exhibit a nonzero duality gap, the coupling between sensitivities and dual variables is invalidated. 
To address this issue, a \emph{convex relaxation} of the SP in~\eqref{SubProb} is proposed in this work to establish a \emph{zero duality gap} in primal-dual optimal solutions for SPs and thus, accurate BFCs. Moreover, the convex relaxation not only ensure convergence, but also improves the computational properties of SPs in terms of scalability and global optimal solvability~\cite{LowConvex1, LowConvex2, Hotz1, Hotz2, LavaeiConvex}.

To this end, let $ X= (x_{i,j})_{N \times N} \in \mathbb{S}^N $. 
Suppose AC line $ k $ connects bus~$ i $ to bus~$ j$ and let $\tilde{x}_{k} = \sqrt{2}x_{i,j}\ \left(  \textrm{hence,}\ (\tilde{x}_{k})^* = \sqrt{2}x_{j,i}\right)$. Therewith, the vectorization $ \xi : \mathbb{S}^N \rightarrow \mathbb{R}^{N+2K} $ of this (partial) Hermitian matrix is defined as
\begin{align}
\xi(X) = [\hspace{0.6mm} &x_{1,1}, \ldots, x_{N,N},\ \real(\tilde{x}_{1}),\dots, \real(\tilde{x}_{K}),\nonumber \\
&\imag(\tilde{x}_{1}), \dots, \imag(\tilde{x}_{K})\hspace{0.6mm} ]^{\mathrm{T} } \in \mathbb{R}^{N+2K}\,. \label{vec}
\end{align}
Let $ W^{t,\omega}$ substitute $ v^{t,\omega}(v^{t,\omega})^{\mathrm{H}} $ and let $\bar{v}^{t,\omega} = \xi(W^{t,\omega})$. Then, the SOC relaxation of~\eqref{ACPF1} can be formulated as in~\eqref{OPF_SOC}, cf. \cite{Hotz2,MohashaPSCC}.
\begin{subequations}\label{OPF_SOC}
	\begin{align}
	& P^{t,\omega}_{\mathrm{G},n} - P^{t,\omega}_{\mathrm{D},n} +P^{t,\omega}_{\mathrm{W},n} \nonumber \\
	&\hspace{1.5cm} = \xi({P}_n^{\mathrm{T}})^{\mathrm{T}} \bar{v}^{t,\omega} + {h}_n^{\mathrm{T}} p^{t,\omega} ;\,   \forall n\in{\mathcal{N}} \label{OPF_SOC_a} \\
	& Q^{t,\omega}_{\mathrm{G},n}- Q^{t,\omega}_{\mathrm{D},n} +Q^{t,\omega}_{\mathrm{W},n} = \xi({Q}_{n}^{\mathrm{T}})^{\mathrm{T}} \bar{v}^{t,\omega};\, \forall n\in{\mathcal{N}} \label{OPF_SOC_b} \\
	& (v_n^{\min})^2 \leq \xi( {M}_n^{\mathrm{T}})^{\mathrm{T}} \bar{v}^{t,\omega} \leq (v_n^{\max})^2 ;\,  \forall n\in{\mathcal{N}} \\
	& \xi( {\hat{I}}_{k}^{\mathrm{T}})^{\mathrm{T}} \bar{v}^{t,\omega} \leq (I_{k}^{\max})^2 ;\, \forall k \in{\mathcal{K}}\\
	& \xi( {\check{I}}_{k}^{\mathrm{T}})^{\mathrm{T}} \bar{v}^{t,\omega} \leq (I_{k}^{\max})^2;\, \forall k \in \mathcal{K}\\
	& S_k(\bar{v}^{t,\omega}) \succeq 0 ;\, \forall k \in{\mathcal{K}} \label{PSD}\\
	& p_l^{\min} \leq p^{t,\omega}_{l} \leq p_l^{\max} ;\,   \forall l \in\mathcal{L} \label{DClim}\\
	& -{r}^{t,\omega}_{{\rm d},n} \leq P^{t,\omega}_{\mathrm{G},n} - {P}^{t}_{\mathrm{G},n} \leq {r}^{t,\omega}_{{\rm u},n};\,  \forall n\in{\mathcal{N}} \label{RES_lim}\\
	& {u}^t_n Q^{\min}_{\mathrm{G},n} \leq Q^{t,\omega}_{\mathrm{G},n} \leq {u}^t_n Q^{\max}_{\mathrm{G},n} ;\,  \forall n\in{\mathcal{N}} \label{Qlim} 
	\end{align}
\end{subequations}
A necessary condition for exactness of the SOC relaxation~\eqref{OPF_SOC} is that all $ 2\times2 $ principal submatrices of  $ W^{t,\omega} $ are \emph{positive semidefinite} (PSD) \cite{Hotz2}. The PSD constraint on the $2\times 2$ principal submatrix $ S_k(\bar{v}^{t,\omega}) $ of $ W^{t,\omega} $ related to AC line~$k$ connecting bus~$ i $ to bus~$ j $ is given by~\eqref{PSD}, which is implemented as an SOC constraint~\cite{Hotz2,MohashaPSCC}. 
The convex feasibility verification SP for scenario~$ \omega \in \Omega$ at time $ t \in \mathcal{T} $ along the lines of \eqref{SubProb} and \eqref{OPF_SOC} can be formulated as in~\eqref{SP_SOC}.
\begin{subequations}\label{SP_SOC}
	\begin{align}
	z^{t,\omega}=\min_{\substack{P^{t,\omega}_{\mathrm{G},n},\, Q^{t,\omega}_{\mathrm{G},n}\\ \bar{v}^{t,\omega},\,{p}^{t,\omega}}}\, & 
	\sum_{n\in{\mathcal{N}}} \sigma^{t}_n\hspace{-1mm}\left( P^{t,\omega}_{\mathrm{G},n},\, {P}^{t}_{\mathrm{G},n} +{r}^{t,\omega}_{{\rm u},n}, \, {P}^{t}_{\mathrm{G},n} -{r}^{t,\omega}_{{\rm d},n}\right)  \nonumber \\
	& \hspace{-3.5mm}+ \sum_{n\in{\mathcal{N}}} \sigma^{t}_n\hspace{-1mm}\left( Q^{t,\omega}_{\mathrm{G},n},\,{u}^t_n Q^{\rm max}_{\mathrm{G},n},\, {u}^t_nQ^{\rm min}_{\mathrm{G},n}\right) \label{SP_obj_SOC} \\
	{\rm s.t.}\ & \eqref{OPF_SOC_a}{\text -}\eqref{DClim}
	\end{align}
\end{subequations}%

\subsection{Exactness and the Hybrid AC/DC Grid Architecture} \label{HTG}
The SOC relaxation \eqref{SP_SOC} of (\ref{SubProb}) is \emph{exact} if and only if the partial matrix $(W^{t,\omega})^\star $ associated with the solution $ (\bar{v}^{t,\omega})^\star $ of~\eqref{SP_SOC} permits a PSD {rank-$ 1 $} completion. For conventional transmission grids, this is typically not the case \cite{Hotz2}. 
The hybrid architecture proposed in \cite{Hotz1}, which comprises a radial AC subgrid with additional DC lines, supports the exactness of the SOC relaxation under normal operating conditions, unless the dual variables of power injections combine to a point in a union of linear subspaces as proven in \cite[Sec. VII]{Hotz2}. 
%
Therefore, the hybrid architecture establishes the accuracy of sensitivities and thereby validates the BFCs. In addition, it ensures the global optimal solvability of SPs~\eqref{SP_SOC} and the applicability of the optimal day-ahead schedule.. 
It shall be noted that the optimal bus voltage vector $ ({v}^{t,\omega})^\star $ associated with the solution $ (\bar{v}^{t,\omega})^\star $ can be recovered using the tree traversal method given in \cite[Sec.~{III-B-3}]{Subhonmesh} and the exactness can be verified a posteriori via the reconstruction error given in~\cite[Sec.~VI-C]{Hotz2}. 
\subsection{Interpretation and Handling of Inexactness}\label{minPL}
As discussed in Section~\ref{HTG}, exactness of the SOC relaxation~\eqref{SP_SOC} is guaranteed if the dual variables of power injections do not form a pathological profile, where the latter are very unlikely in case that the cost function has a positive gradient of the optimizer~\cite{Hotz2}. This is mostly the case in initial iterations owing to the cost of up-reserve allocation. However, the aforementioned condition is not satisfied when the day-ahead schedule is feasible (then, the solution lies in the second segment of~\eqref{SP_obj_SOC} which has a zero gradient) or when the down-reserves are insufficient (then, the solution  lies in the third segment of~\eqref{SP_obj_SOC} which has a negative gradient).
In this regard, a convex auxiliary problem (AP) is defined as in~\eqref{SP_Feas} to verify the feasibility of the day-ahead schedule for scenario~$ \omega $ at time~$ t $.
\begin{align}\label{SP_Feas}
\min_{\substack{P^{t,\omega}_{\mathrm{G},n},\, Q^{t,\omega}_{\mathrm{G},n}\\ \bar{v}^{t,\omega},\,{p}^{t,\omega}}}\,  \Bigg\{
\sum_{n\in{\mathcal{N}}} P^{t,\omega}_{\mathrm{G},n} \
\Bigg\arrowvert \	\eqref{OPF_SOC} \Bigg\}
\end{align}
Therein, the objective function minimizes the power loss  subject to the complete set of AC-NCs for scenario $ \omega $ at time $ t $. Therefore, the feasibility of~\eqref{SP_Feas} confirms the validity of the day-ahead schedule generated by~\eqref{UC_Prob} for scenario~$ \omega $ at time~$ t $, \emph{i.e.}, $ z^{t,\omega} = 0 $ and no BFC is required.

On the other hand, if \eqref{SP_Feas} is infeasible, that reflects the insufficiency of down-reserve allocation. 
This is heuristically rectified by re-evaluating the MP with an adjusted $ \alpha^{t} \leftarrow \Upsilon\alpha^{t} $ in~\eqref{dnRes}, using a predetermined constant~$ \Upsilon>1 $.
%
\subsection{Summary of the Proposed Algorithm}
Finally, Algorithm~\ref{Algo} summarizes the solution strategy to compute the robust day-ahead optimal schedule $ \Phi^\mathcal{T}_\mathcal{N} $.
\begin{algorithm}[h!] 
	\SetKwInOut{Init}{Initialize}
	\Init{$\eta=0,\, \Upsilon=1.1,\, \epsilon=0.005,\,\alpha^t=1;\, \forall t \in \mathcal{T}$ and $ P^{t,\omega}_{{\rm loss}}=0,\, z^{t,\omega}=\infty;\,\forall \omega \in \Omega,\, \forall t \in \mathcal{T}$.}
	\While{$ \underset{t \in \mathcal{T},\omega \in \Omega}{\max}z^{t,\omega} \geq \gamma \epsilon$}{
		$ \eta\leftarrow \eta+1 $.
		\\
		Execute the MP~\eqref{UC_Prob} and extract $ \Phi^\mathcal{T}_\mathcal{N}$.
		\label{S1}\\
		\For{$ t \in \mathcal{T} $} 
		{
			\For{$ \omega \in \Omega $}{			
				Execute the SP~\eqref{SP_SOC}.
				\\
				\eIf{Exact}{Calculate $ P^{t,\omega}_{\rm loss} $ using \eqref{P_LOSS}.}{Execute the AP~\eqref{SP_Feas}.\\
					\eIf{Feasible}{Calculate $ P^{t,\omega}_{\rm loss} $ using \eqref{P_LOSS}.}{Adjust $ \alpha^{t} \leftarrow \Upsilon\alpha^{t} $. Return to Line~\ref{S1}.} }
			}
			Compute the BFC as in~\eqref{BenCut}.
		} 
		Update $ P^{t,\omega}_{\rm loss},\,  \forall \omega \in \Omega,\, \forall t \in \mathcal{T}  $.\\
		Add BFCs into the MP.
	}
	\caption{Robust day-ahead optimal energy and reserve management.} \label{Algo}
\end{algorithm}\vspace{-0.5em}%

\section{Simulation Results} \label{Results}
\subsection{Test System} \label{system}
The proposed day-ahead energy and reserve optimization method is illustrated using the $ 2383 $-bus test case “case$ 2383 $wp.m” which represents the Polish transmission grid during winter peak conditions. The test case is provided by the power system simulation package MATPOWER~\cite{MATPOWER} and is preprocessed as in \cite[Sec.~VIII-A]{Hotz2} to support the simulations. Further, the details about the upgrade strategy to a hybrid transmission grid (HTG) can be found in \cite[Sec.~VIII-B]{Hotz2}. The day-ahead energy and reserve optimization of the original AC transmission grid (ACG) and the HTG are compared in simulations. In addition to the data of the test case, the following parameters related to the UC formulation~\eqref{UC_Prob} are assumed:  $ P^{\min}_{{\rm G},n}={\max}\{P^{\min}_{{\rm G},n},10\,{\rm MW}\} $, $ c^{\rm e}_{0,n}= \$20/{\rm h}$, $ c^{\rm s}_n= \$100$, $ c^{\rm d}_n= \$10$, $ R^\Delta_{{\rm u},n}=R^\Delta_{{\rm d},n}=0.5(P^{\max}_{{\rm G},n}-P^{\min}_{{\rm G},n}) $, $ \hat{T}_n=4\,{\rm h} $, $ \check{T}_n=2\,{\rm h} $ and $ R^\delta_{{\rm u},n}=R^\delta_{{\rm d},n}=0.25(P^{\max}_{{\rm G},n}-P^{\min}_{{\rm G},n}) $, for all $ n \in \mathcal{N} $. 
The day load profile which is used to scale the individual active power demand at each bus over the time horizon, is taken from \cite[Table~IV]{RTS}. Similarly, the reactive power demand is also scaled assuming a constant power factor over the entire time horizon for all loads. 
Further, the wind profiles of $ 15 $ wind farms from~\cite{RES_Wind_2016} are integrated at buses $ 6 $, $ 8 $, $ 9 $, $ 15 $, $ 31 $, $ 32 $, $ 183 $, $ 682 $, $ 711 $, $ 723 $, $ 729 $, $ 833 $, $ 1230 $, $ 1283 $ and $ 1546 $. The peak value of all wind profiles is adjusted to $ 140\,{\rm MW} $ and the rest of the values are scaled proportionally. The load buses are grouped into $ 15 $ clusters where the demand uncertainty within each cluster is assumed to be equal and the variations are associated to a column in the selected OA. 
To this end, the first $ 30 $ columns of the OA $ L_{32}2^{31} $ are adopted to form $ 32 $ scenarios except the base-case for each time instant~$ t \in \mathcal{T} $~\cite{TOAT}. Each scenario represents information on either of the confidence bounds (as defined in~\eqref{Pd}) of the $ 30 $ uncertain variables of the respective time instant~$ t $.
Therein, $ P^{t,\zeta-}_{{\rm D},n}=P^{t,\zeta+}_{{\rm D},n}=0.05{P}^{t}_{{\rm D},n}$, $ P^{t,\zeta-}_{{\rm W},n}=0.5{P}^{t}_{{\rm W},n} $ and $ P^{t,\zeta+}_{{\rm W},n}=0.1{P}^{t}_{{\rm W},n} $, for all $ n \in \mathcal{N} $ in all case studies.%
\vspace{-0.5em}
\subsection{Numerical Simulations}\label{Simulations} 
The MP in~\eqref{UC_Prob} is a mixed-integer linear program which is solved using CPLEX. For the HTG, the SP in~\eqref{SubProb} are convexified using the SOC relaxation discussed in Section~\ref{ConvRel}. The convex SPs~\eqref{SP_SOC} and APs~\eqref{SP_Feas} are solved using the primal-dual interior-point solver MOSEK. However, \eqref{SP_SOC} is inexact for the ACG as discussed in Section~\ref{HTG}~\cite{Hotz2}. Consequently, its network constraints and hence~\eqref{SubProb} is approximated with the decoupled network representation for the ACG. The resulting linear program is solved using CPLEX. However, AC power flow simulations are performed on the optimal schedule to verify the viability of the base-case scenario at every hour. The base value used in computations is $ 100\,{\rm MVA} $.%
\subsection{Case Study 1: Economic Efficiency}
\begin{figure}[th]
	\centering
	\includegraphics[width=3.2in]{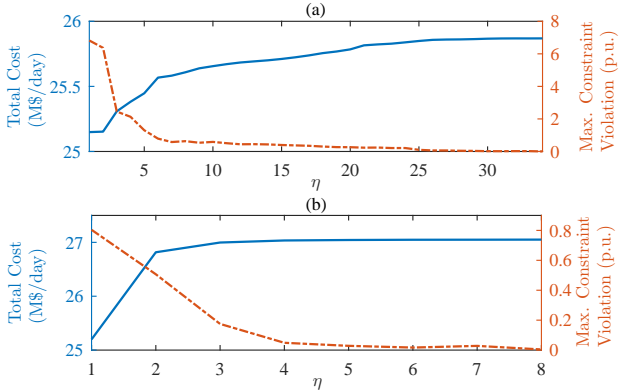}\vspace{-0.5em}%
	\caption{Convergence trajectories of the algorithm for (a) the ACG and (b) the HTG in terms of the objective value and constraint violations.}%
	\label{LB_Cov}
\end{figure}

The following simulation results verify that the iterative algorithm performs satisfactorily for both grids. Fig.~\ref{LB_Cov} depicts the convergence trends at the normal demand. The day-ahead schedule is corrected during the iterations to avoid constraint violations. Consequently, the total cost of the base-case operation (objective value of the MP) increases in successive iterations. The algorithm converges faster for the HTG as its flexibility induces less constraint violations and improves the accessibility of economic generation facilities.%

\begin{figure}[th]
	\centering
	\includegraphics[width=3.2in]{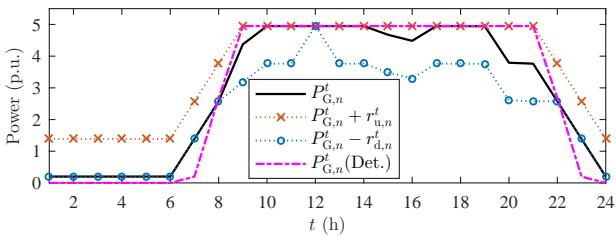}\vspace{-0.5em}%
	\caption{Schedule for generator $ 175 $ in the HTG. Curve~$ 1 $ is the base-case power dispatch. Curve~$ 2 $ and $ 3 $ are the up and down reserve margins for the robust operation. Curve~$ 4 $ is the power dispatch for a deterministic operation.}\label{Pg} 
\end{figure}

In the ACG and HTG, $ 314 $ and $ 315 $ units are committed with $ 259 $ and $ 263 $ units constantly operating at their maximum capability respectively.
For instance, Fig.~\ref{Pg} illustrates the schedule for generator $ 175 $ in the HTG.
It can be observed that for a deterministic operation, it is committed only from hour~$ 7 $ to hour~$ 23 $. In contrast, it is committed over the entire time horizon for the robust operation, being partially dispatched at aforementioned hours to support reserve requirements.%

Table~\ref{GridLoad} reports the total generation cost for the deterministic and robust operation at different system loads 
for both grids. Furthermore, it compares the cost reduction for a robust operation in the HTG with respect to the ACG. The total cost for the ACG is calculated including the generation cost of the slack power. It can be observed that the HTG offers a better utilization of generation resources, thereby enabling consistent economic benefits at different demands in comparison to the ACG. Although the generation capability is sufficient, the ACG \emph{cannot} be loaded further due to transmission congestion issues. In contrast, the additional flexibility of the HTG in power flow routing results in improved loadability up to $ 14.5\% $ of the normal demand, \emph{i.e.}, it can accommodate increases in demand.
%
%
\begin{table}[!h]
	\renewcommand{\arraystretch}{1.3}
	\centering
	\caption{Comparison of Economic Efficiency w.r.t. Grid Loadability}
	\label{GridLoad}
	\begin{tabular}{c|c c|c c|c}
		\hline
		Scale	& \multicolumn{2}{c|}{Deterministic}&\multicolumn{3}{c}{Robust}\\ \cline{2-6}
		Factor & \multicolumn{4}{c|}{Total Cost $ ({\rm M }\$/\mathrm{day}) $} & Cost \\  \cline{2-5}
		& ACG & HTG & ACG & HTG & Reduction\\
		\hline \hline
		$  0.95 $ & $ 23.79 $ & $ 23.57 $ & $ 24.41 $ & $ 23.94 $ & $ 1.94\% $\\
		$ 1.00 $ & $ 26.91 $ & $ 26.65 $ & $ 27.75 $ & $ 27.05 $ & $ 2.51\% $\\
		$ 1.02 $ & $ 28.23 $ & $ 27.92 $ & $ 29.54 $ & $ 28.37 $ & $ 3.96\% $\\
		$ 1.145 $ & $ - $ & $ 36.24 $ & $ - $ & $ 36.87 $ & $ - $\\
		\hline
	\end{tabular}\vspace{-1.0em}
\end{table}
Further, it should be noted that the generation schedules of the ACG, which are computed based on a decoupled representation of the network, \emph{cannot} be implemented due to a significant number of violations in voltage limits, transmission line flow limits, and reactive power capabilities of generators according to AC power flow simulations performed with MATPOWER~\cite{MATPOWER}. Consequently, the operation cost of the ACG are potentially even higher than the values listed in Table~\ref{GridLoad}.

\subsection{Case Study 2: Generation Utilization}
In this case study, the wind penetration is increased by $ 25\% $ and the generation utilization is examined against grid loadability. The simulation of the ACG reported \emph{infeasibility}, indicating that the ACG \emph{cannot} accommodate this much wind-based injections. This reflects the inflexibility of the ACG to deploy reserves owing to the transmission congestion, despite the fact that the available generation capacity is sufficient. The results for the HTG are documented in Table~\ref{GenUti}. The total cost is increasing at an almost linear rate and all the generators are committed for $ 1.05 $ or more times the normal demand. The energy and reserve contributions are computed as a percentage of the generation capacity of online units, \emph{i.e.}, $ \underset{t \in \mathcal{T}}{\sum}\underset{n \in \mathcal{N}}{\sum} u^t_nP^{\max}_{{\rm G},n} $. Around $ 70\% $ to $ 80\% $ of online capacity is used to meet the forcasted energy demand and the utilization for reserves is less than $ 10\% $. The main reason for the partial utilization of the online capacity is the ramping limit of the units.%
\begin{table}[!h]
	\renewcommand{\arraystretch}{1.3}
	\centering
	\caption{Generation Utilization of the HTG w.r.t. Grid Loadability at High Wind Power Penetration}
	\label{Pw}\label{GenUti}
	\begin{tabular}{c c c|c c c}
		\hline
	Scale & Total Cost & \# Online & \multicolumn{3}{c}{Generation Utilization}\\ \cline{4-6}
	Factor & $ ({\rm M }\$/\mathrm{day}) $ & Units & Energy & Up-Res. & Down-Res.\\  
		\hline \hline
		$  0.95 $ & $ 22.91 $ & $ 311 $ & $ 72.53\% $ & $ 8.12\% $ & $ 6.51\% $\\
		$ 1.00 $ & $ 26.00 $ & $ 316 $& $ 74.70\% $ & $ 8.06\% $ & $ 6.02\% $\\
		$ 1.05 $ & $ 29.27 $ & $ 321 $& $ 76.20\% $ & $ 7.93\% $ & $ 6.39\% $\\
		$ 1.10 $ & $ 32.54 $ & $ 321 $& $ 77.83\% $ & $ 7.98\% $ & $ 7.12\% $\\
		$ 1.145 $ & $ 35.71 $ & $ 321 $& $ 80.30\% $ & $ 8.11\% $ & $ 8.61\% $\\
		\hline
	\end{tabular}
\end{table}\vspace{-0.1em}%

\begin{figure}[th]
	\centering
	\includegraphics[width=3.2in]{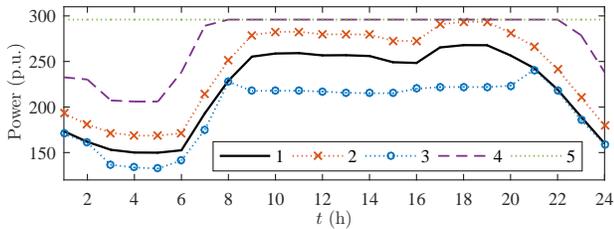}\vspace{-0.9em}%
	\caption{Aggregated generation profiles of the HTG at $ 14.5\% $ increased demand. Curve~$ 1 $ is the total power dispatch for the base-case. Curve~$ 2 $ and $ 3 $ are the up and down reserve margins. Curve~$ 4 $ is the generation capacity of the online units. Curve~$ 5 $ is the total generation capacity of the system.}%
	\label{SumP} \vspace{-0.5em}
\end{figure}

Fig.~\ref{SumP} illustrates the aggregated generation profile of the system at $ 114.5\% $ of the normal demand. 
It can be observed that at hour~$ 18 $ (peak hour), the total generation capacity is almost exploited (up to $99.19\% $) for energy and up reserve requirements. Hence, the flexibility offered by the HTG supports the effective utilization of \emph{all} generation in the system in contrast to ACG. 
%
%
\vspace{-0.3em}%

\section{Conclusion} \label{Conclu}
This paper discussed the electrical system model and mathematical foundation of the UC problem with AC-NCs under uncertainty. Therein, confidence bounds are used to characterize the  demand and RES uncertainty and orthogonal arrays are used to capture the worst-case uncertain scenarios. A two-stage iterative algorithm was proposed based on the BDF in which the day-ahead optimal schedule is computed in the first stage and the second stage verifies its viability against the AC-NCs for possible uncertain realizations. In addition, a convex relaxation was applied to the SPs which enables their globally optimal solution and supports convergence. Numerical simulations were performed for the Polish $ 2383 $-bus system during winter peak conditions. Therein, the merits of the recently proposed hybrid transmission grid architecture were illustrated in comparison to the original grid in a day-ahead market context. Firstly, it ensures exactness of the convex relaxation of SPs, thereby validating the optimality and the applicability of the day-ahead schedule. Secondly, the flexibility offered by the architecture alleviates transmission congestion issues, improving the utilization of the generation resources. This results in lower total cost over different load levels. Further, its enhanced loadability and reserve provision enable the utilization of \emph{all} generation in the system and to accommodate the increase in demand and RESs. In that respect, the hybrid architecture can be utilized as a topology-preserving capacity expansion strategy for existing congested grids, which supports an optimal and complete utilization of available generation facilities.
\bibliographystyle{IEEEtr}
\bibliography{references_HTG} 
%

%






\end{document}